\documentclass[12pt]{article}
\input xy
\xyoption{all}

\usepackage{amsmath}
\usepackage{amssymb}
\usepackage{amsthm}

\begin{document}
\input{amssym.def}

\newsymbol \circledarrowleft 1309

\newtheorem{guess}{Theorem}
\newcommand{\bth}{\begin{guess}$\!\!\!${\bf }~}
\newcommand{\eeth}{\end{guess}}

\newtheorem{lema}[guess]{Lemma}
\newcommand{\blem}{\begin{lema}$\!\!\!${\bf }~}
\newcommand{\elem}{\end{lema}}

\newtheorem{Cor}[guess]{Corollary}
\newcommand{\bcor}{\begin{Cor}$\!\!\!${\bf }~}
\newcommand{\ecor}{\end{Cor}}

\newtheorem{Def}[guess]{Definition}
\newcommand{\bdefe}{\begin{Def}$\!\!\!${\bf }~}
\newcommand{\edefe}{\end{Def}}

\newtheorem{Prop}[guess]{Proposition}
\newcommand{\bprop}{\begin{Prop}$\!\!\!${\bf }~}
\newcommand{\eprop}{\end{Prop}}

\newtheorem{eg}[guess]{Example}
\newcommand{\beg}{\begin{eg}$\!\!\!${\bf }~\rm}
\newcommand{\eeg}{\end{eg}}

\newtheorem{rema}[guess]{\bf Remark}
\newcommand{\brem}{\begin{rema}$\!\!\!${\bf }~\rm}
\newcommand{\erem}{\end{rema}}

\newcommand{\ra}{\rightarrow}
\newcommand{\lr}{\longrightarrow}
\newcommand{\ct}{{\cal T}}
\newcommand{\cp}{{\cal P}}
\newcommand{\cz}{{\cal Z}}
\newcommand{\cb}{{\cal B}}
\newcommand{\cy}{{\cal Y}}
\newcommand{\cf}{{\cal F}}
\newcommand{\cl}{{\cal L}}
\newcommand{\cx}{{\cal X}}
\newcommand{\bc}{{\Bbb C}}
\newcommand{\bp}{{\Bbb P}}
\newcommand{\br}{{\Bbb R}}
\newcommand{\ba}{{\Bbb A}}
\newcommand{\bz}{{\Bbb Z}}
\newcommand{\bg}{{\Bbb G}}

\newcommand{\ctext}[1]{\makebox(0,0){#1}}
\setlength{\unitlength}{0.1mm}

\newcommand{\ol}{\overline}
\newcommand{\spec}{{\rm Spec}\,}
\newcommand{\wt}{\widetilde}

\title{Toric degeneration of Bott-Samelson-Demazure-Hansen varieties}
\author{\footnote {Paramasamy Karuppuchamy, Department of Mathematics
and Statistics, Mail Stop 942,The University of Toledo
2801 W. Bancroft St.,Toledo, OH 43606-3390, USA.
  Email: paramasamy@gmail.com} Paramasamy Karuppuchamy and
\footnote {School of Mathematics, Tata Institute of Fundamental
Research, Homi Bhabha Road, Mumbai 400005, India. Email:
param@math.tifr.res.in}
  Parameswaran, A. J.}
\date{}
\maketitle

{\abstract In this paper we construct a degeneration of
Bott-Samelson-Demazure-Hansen varieties to toric varieties in an
algebraic family and study the geometry of the resulting toric
varieties. We give a natural set of torus invariant curves that
generate the Chow group of $1$-cycles of the limiting toric variety
and express the ample cone of this toric variety as a sub-cone of
the ample cone of the corresponding Bott-Samelson-Demazure-Hansen
variety. We also give a description of Extremal and Mori rays and determine 
when this toric variety is Fano. }

\section{Introduction}

Let $ G$ be an almost simple, simply connected, affine algebraic group
defined over an algebraically closed field $ k$ of arbitrary
characteristic. Let $B$ be a Borel subgroup of $G$.
Then $G$ acts from the left on the flag variety $G/B$. The
$B$-invariant closed subvarieties of $G/B$ are called {\it Schubert
varieties}.  Every Schubert variety is uniquely represented by an
element of the Weyl group. After choosing a reduced expression
for  Weyl group elements as product of simple
reflections, one constructs certain smooth birational modifications of
the corresponding Schubert varieties. These desingularizations used to be
 called as Bott-Samelson varieties and the constructions were first described by Demezure and 
indepently by Hansen (for more details \cite{jantzen}, Chapetr 13, page 353). These constructions
 (cf. \cite{dem} or \cite{jantzen})  extend
naturally to any given sequence of simple reflections. The resulting varieties are referred as
 Bott-Samelson-Demazure-Hansen varieties in this paper.
We often abbreviate this long name and call these as BSDH varieties.

For $k=\bc$, the field of complex numbers, the authors Grossberg and
Karshon (cf. \cite{gross}) obtained a family of complex structures on
BSDH variety as differentiable manifolds which
degenerate to a different complex structure, with resulting manifold
having a toric variety structure. We give an algebraic degeneration by
constructing a smooth family of varieties parametrised by the affine
line with general fibre isomorphic to the
BSDH variety and the special fibre isomorphic
to a smooth toric variety (cf. Section $3$).

In another direction it may be desirable to obtain degenerations of
Schubert varieties using the degenerations of their
BSDH resolutions. Using Standard Monomial
basis Gonciulea and Lakshmibai \cite{goncu} degenerated Flag varieties
and some Schubert varieties to toric varieties.  Later Caldero
\cite{cald} had constructed degeneration of Schubert varieties into
toric varieties using Lusztig's canonical basis.

Lauritzen and Thomsen  have given an ampleness
criterion for line bundles on these BSDH
varieties (cf. \cite{laurit}, Theorem 3.1). In fact,  they gave a set of line bundles and showed
that the ample cone is the strict
positive cone generated by these line bundles ({cf. section~\ref{Amplecone}}).
Here we describe the ample cone of the toric
variety.   We also give a necessary and sufficient condition for the anti canonical
bundle of the limiting toric variety to be ample, i.e. we give criterion to decide whether the limiting toric variety is a Fano variety.

The paper is arranged as follows. In section $2$ we fix some
notations and recall some basic results. The main observation here
is the relation between the self intersection number of a
BSDH surface constructed using two simple
reflections and the pairing between the corresponding roots
(Lemma~\ref{basic}). Section $3$ is devoted to obtain a
degeneration of BSDH variety to a toric
variety in an algebraic family (Theorem \ref{piprops}). In section
$4$ we study the Chow group of cycles of dimension $1$ (curves) of
the limiting toric variety. We label certain torus invariant curves
and obtain some basic relation between them
(Proposition~\ref{formula}, Remark~\ref{H-surfaces}). In section
$5$, we obtain a basis  for the Chow group of one cycles of the limiting toric variety such that 
all torus invariant curves are non negative linear combination of elelments of this basis (Theorem
\ref{invariantlines}). We also give two algorithms to find these
torus invariant curves (Lemmas~\ref{compalgo} and~\ref{weylalgo}).
In section $6$ we study the Extremal rays and Mori rays and
describe them completely for the limiting toric variety
(Theorem~\ref{Extrmalrays}, Theorem~\ref{mori}). We also give a criterion to decide whether the limiting toric variety is Fano or not.
The last section deals with description of the ample cone of these toric varieties
as a subcone of the ample cone of the corresponding
BSDH variety (Theorem \ref{amplecone}).

\section {Preliminaries}

Let $ G$ be an almost simple, simply connected, affine algebraic
group defined over an algebraically closed field $ k$ of arbitrary
characteristic.  Fix a maximal torus $T\subset G$.  Then the Weyl
group is defined as $N(T)/T$, where $N(T)$ is the normaliser of $T$
in $G$. If we denote the character group of $T$ by $X(T)$, then the
Weyl group $W$ has a faithful representation on the real vector
space $X(T)\otimes\br$. Let $(~,~)$ be a non-degenerate $W$ invariant pairing on $X(T)\otimes\br$.

Let $ \ S, \ \Phi^+, \ \Phi^- \subset X(T)$ be simple roots, positive
roots and negative roots respectively.  Let $U_\alpha$ be the root subgroup corresponding to the root
$\alpha$ and $B$ be the Borel subgroup
of $G$ generated by the root subgroups corresponding to the negative
roots and maximal torus $T.$ A closed subgroup of $G$ containing a Borel subgroup is called a parabolic subgroup.  
 For a simple root $\alpha$, the minimal parabolic subgroup
$P_\alpha$ is defined to be the subgroup generated by $B$ and the root
subgroup $U_\alpha$. The fundamental weight corresponding to the simple root $\alpha$ is 
denoted by $\omega_\alpha$.

Let $\alpha^\vee:=\frac{2\alpha}{(\alpha,\alpha)}$ be the co-root of
$\alpha$.  For a given simple root $\alpha \in S$, let $s_{\alpha}$
denote the simple reflection on $X(T)\otimes\br$ defined by
$s_{\alpha}(x) = x-(x,\alpha^\vee)\alpha$. Then $W$ is generated by
$\{s_\alpha \mid \alpha\in S\}$.  

Recall the following well-known description of parabolic subgroups of
$G$. Let $\xi :{\bg}_ m \lr G$ be a one parameter subgroup of $G$. 
% Define a ${\bg}_m$-action on $G$ by $ x \ast g = \xi (x) g \xi (x)^{-1}.$ 
Then we have the following (cf. \cite{spring} p. 148):

\blem \label{parab} The set $ P(\xi):=\{g \in G~~|~ \lim_{x
\rightarrow 0} \ \xi (x) g \xi (x)^{-1}\}$ is a parabolic subgroup
and the unipotent radical $ R_u(P(\xi))$ of $ P(\xi)$ is given by\\
$\{g \in G~~|~ \lim_{x \rightarrow 0}\ \xi (x) g \xi (x)^{-1} = {\rm identity}\}$.

Moreover any parabolic subgroup of $G$ is of the form $P(\xi)$ for
some $\xi :{\bg}_ m \lr G$.
\hfill$\diamondsuit$ \elem

We fix a $\xi:{\bg}_ m \lr T \subset G$ once and for all such that $B = P(\xi)$
and the unipotent radical $B_U$ of $B$ is $ \{ g \in G ~|~ \lim_{x
\rightarrow 0} \xi (x) g \xi (x)^{-1} =$ identity $\}$. Then we have the following:

\blem\label{familyhom} Let $\ba:=spec~k[t]$ denote the affine line
over $k$. Define a family of homomorphisms $\phi_x:B\to B$
parametrised by $x\in \ba$ as follows: $$ \phi_x(b) = \xi
(x) b \xi (x)^{-1} ~~{\rm if} ~x\neq 0 ~~~~~ {\rm and}~~~~ \phi_0(b) =
\lim_{x \rightarrow 0}~ \xi (x) b \xi (x)^{ - 1} $$ Then $\phi_x$ is
an automorphsm for $x\neq 0$ and $\phi_0$ is the natural Levi
projection $B\to T$ with unipotent radical $B_U$ of $B$ as the kernel.
\hfill$\diamondsuit$ \elem

Let ${\bf w} = (s_1, s_2, \cdots , s_m)$ be a sequence of simple
reflections, where $m$ is any positive integer.  Each simple reflection $s_j$ is defined by a simple 
root, say $\alpha _j$. We call $m$ the {\em length of the sequence} and we denote it by  $l({\bf w})$. Length
 of an empty sequence is defined to be $0$. Note that this lenth is just the number of terms in the 
sequence.  We also have a lenth function on the Weyl group $W$. Every element $w$ of $W$ can be 
written as a product of simple reflections, $w= s_{j_1}s_{j_2}\cdots s_{j_r}$. If $w$ can not be 
written as of less than $r$ number of simple reflections then this expression is called a {\em reduced 
expression} and $r$ is called the the length of $w$. We have a natural map from (finite) sequences of 
simple relfections to the the Weyl group $W$.  This maps a sequence  
${\bf w} = (s_1, s_2, \cdots , s_m)$ to the element $w=s_1s_2 \cdots s_m$. But this representation 
may not be reduced and hence the length of the sequnce ${\bf w}$ 
may not be equal to the length  of the corresponding Weyl group element $w$.

 A sequence of integers $I=(i_1,\cdots, i_r)$ is called
{\em $m$-admissible} if $1\leq i_1<i_2<\cdots<i_r\leq m$. The entries
$i_1$ and $i_r$ of $I=(i_1,\cdots, i_r)$ are called {\em initial} and
{\em final} entries respectively. Define the subsequence ${\bf w}_I$
of  ${\bf w} = (s_1, s_2, \cdots , s_m)$ for every $m$-admissible sequence $I=(i_1,\cdots, i_r)$
by ${\bf w}_I:=(s_{i_1},\cdots, s_{i_r})$. For $0 \leq r \leq m$ we define  the truncated $m$-admissible sequences $I[r]:=(1,~ 2,~\cdots ~, r)$ and
$[r]I:=(r+1,r+2,\cdots,m)$, we denote the corresponding subsequence of
simple reflections ${\bf w}_{I[r]}$ and ${\bf w}_{[r]I}$ by ${\bf
w}[r]$ and $[r]{\bf w}$ respectively.  Note that ${\bf w}[m]={\bf  w}=[0]{\bf w}$ and  $[m]{\bf w}={\bf w}[0]$ is the empty sequence of simple reflections that corresponds to the empty sequences of integers $[m]I=I[0]$.

 Let $P^{\bf w}$ and $B^{\bf w}$ denote the
products $P_1 \times \cdots  \times P_m$ and $B\times \cdots \times
B$ (m copies) respectively, where $P_j$  denotes the minimal
parabolic subgroup $P_{\alpha_j}$.  The Bott-Samelson-Demazure-Hansen (BSDH) variety, $Z_{\bf w}$, is defined ( see
Definition~\ref{bottaction} and Definition~\ref{bsdschemes} below for a
twisted version) as the quotient

$$Z_{\bf w} := \frac {P_1 \times P_2 \times \cdots \times P_m}{B \times B \times \cdot \times B} =\frac{P^{\bf w}}{B^{\bf w}}$$

where the $B^{\bf w}$ acts on $P^{\bf w}$ from the right as follows:  
$$ (p_1,p_2,p_3\cdots p_m)(b_1,b_2,b_3 \cdots b_m)=(p_1b_1, b_1^{-1}p_2b_2, b_2^{-1}p_3b_3 \cdots b_{m-1}^{-1}p_mb_m) $$

The BSDH variety also has the following inductive geometric construction.  
The induction is on the length of the sequence  ${\bf w} = (s_1, s_2, \cdots , s_m)$.  
We construct the BSDH variety $Z_{{\bf w}[r]}$ and a map $f_r: Z_{{\bf w}[r]}\longrightarrow G/B$
inductively for all $0\leq r \leq m$.
 
When $l({\bf w})=0$ i.e. ${\bf w}={\bf w}[0]$, we define the corresponding BSDH variety
to be the unique $B$ fixed point, the identity coset $eB$, of $G/B$. The map $f_0 :Z_{{\bf w}[0]} 
\longrightarrow G/B$ is the inclusion.

When the $l({\bf w})=1$ i.e. ${\bf w}={\bf w}[1]={s_1}$, we define the corresponding BSDH variety 
$Z_{{\bf w}[1]} \cong P_1/B $ as the fiber product $Z_{{\bf w}[0]}\times _{G/P_1} G/B$ 
Observe that our construction gives the BSDH variety $Z_{{\bf w}[1]}$ with the map 
$f_1:Z_{{\bf w}[1]}\longrightarrow 
G/B$, the projection $\psi_1 :Z_{{\bf w}[1]} \longrightarrow Z_{{\bf w}[0]}$, and the section
 $\sigma_0: Z_{{\bf w}[0]} \longrightarrow Z_{{\bf w}[1]}$.

Suppose the BSDH variety $Z_{{\bf w}[r]}$ together with a map $f_r: Z_{{\bf w}[r]}
\longrightarrow G/B$ is already constructed.
  Now the fiber product $Z_{{\bf w}[r]}\times _{G/P_{r+1}} G/B$  defines the BSDH variety 
$Z_{{\bf w}[r+1]}$ and the map  $f_{r+1}$ and $\psi_{r+1}$ are  canonical projections. 

A section $\sigma _r$ to this projection $\psi_{r+1}$ is equivalent to giving a lift
of the map $Z_{{\bf w}[r]}\lr G/P_{r+1}$ to $G/B$.  Our inductive procedure already provided us such a map, 
$f_r$ (see diagram below).

$$ \xymatrix{
 Z_{{\bf w}[r+1]} \ar[rr]^{f_{r+1}} \ar[dd]^{\psi_{r+1}} & & G/B
\ar[dd]^{\pi_{r+1}} \\
& { \square } & \\
Z_{{\bf w}[r]} \ar @/^1pc/ [uu]^{\sigma_r} \ar[rr]^{\pi_{r+1}\circ f_r} & & G/P_{r+1} } $$

In summary, we have indcutively constructed the following:

\noindent
i) The BSDH variety $Z_{{\bf w}[r]}$ with the map $f_r:Z_{{\bf w}[r]}\longrightarrow G/B$, for all $r$. 

\noindent
ii) The projection $\psi_r :Z_{{\bf w}[r]} \longrightarrow Z_{{\bf w}[r-1]}$, for 
$1 \leq r \leq m$. 

\noindent
iii) The section $\sigma_r: Z_{{\bf w}[r]} \longrightarrow Z_{{\bf w}[r+1]}$, for  
$0 \leq r \leq m-1$.

The next two Lemmas appear in [\cite{kempf}, Lemma 3(a), Lemma 2(3)]  in  slightly different 
notations. Since these Lemmas are crucial for our computations we give a proof of one of 
these  lemmas for the convenience of the reader. 

\blem \label{basic} Let ${\bf u}= (s_1, s_2)$ be a sequence
of simple reflections. Then the self intersection number of the section
$\sigma_1 (Z_{s_1})$ in the surface $Z_{\bf u}$ is $(\alpha_2, ~ \alpha_1
^{\vee})$.  By abuse of notation we sometimes denote this number $(\alpha_2, ~ \alpha_1
^{\vee})$ by $(2,1)$.\elem

\noindent {\bf Proof}: The self intersection number of the section
$\sigma _1(Z_{s_1})$ in the surface $Z_{\bf u}$ is by definition the
degree of the normal bundle $N_{\sigma _1(Z_{s_1})/Z_{\bf u}}$. But the
normal bundle of a section in a fibration can be identified with the
restriction of the relative tangent bundle.
Since this fibration is a fibre product of
$Z_{s_1}\rightarrow G/P_2$ with the natural fibration
$\pi_2: G/B\to G/P_2$, the relative tangent bundle is the
pull back of the relative tangent bundle of $\pi_2$

$$
\xymatrix{
Z_{s_2} \ar@{=}[r] \ar[d] & P_2/B \ar[d]\\
 Z_{\bf u} \ar[r]^f \ar[d]^{{\psi}_2}
  & G/B \ar[d]^{\pi_2} \\
   Z_{s_1}\ar @/^1pc/[u]^{\sigma_1}  \ar[r] & G/P_2     }
$$

Now the relative tangent bundle of $\pi_2$ is canonically
identified with $L_{\alpha_2}$, the line bundle on $G/B$ associated to
the character $\alpha_2$.  Hence it suffices to prove that the degree
of $L_{\alpha_2}$ restricted to $f(\sigma_1(Z_{s_1}))$ is $
(\alpha_2, \alpha _1^{\vee})$, as $f$ defines an embedding. In fact
$f(\sigma_1(Z_{s_1}))=P_1/B\subset G/B$.

Using Part II, Proposition 5.2 of [J], we obtain the individual
cohomology groups of the restriction of $L_\lambda$ for any character
$\lambda$ on $P_1/B$, which in turn determine the line bundle.
The Lemma will then follow by substituting the character $\alpha_2$ for
$\lambda$.

Proposition 5.2 (b) of [J] implies when $(\lambda, \alpha_1 ^{\vee})=-1$ then
$H^i(P_1/B,
L_\lambda )=0$ for all $i\geq 0$, hence $L_\lambda\mid_{P_1/B}$ is
isomorphic to ${\cal O}_{P_1/B}(-1)$, whose degree is $-1=(\lambda,
\alpha_1 ^{\vee})$,

Proposition 5.2 (c) of [J] implies when $(\lambda, \alpha_1 ^{\vee})\geq 0$, then
$H^1(P_1/B,
L_\lambda )=0$ and $H^0(P_1/B, L_\lambda )$ has dimension
$(\lambda, \alpha_1 ^{\vee})+1$. Hence the line bundle
$L_\lambda\mid_{P_1/B}$ is
${\cal O}_{P_1/B}((\lambda, \alpha_1 ^{\vee}))$.

Proposition 5.2 (d) of [J] implies  when $(\lambda, \alpha_1 ^{\vee})\leq -2$, then
$H^0(P_1/B, L_\lambda )=0$ and $H^1(P_1/B, L_\lambda )$ has
dimension $-(\lambda,\alpha_1^{\vee})-1$. Hence the line bundle
$L_\lambda\mid_{P_1/B}$ is ${\cal O}_{P_1/B}((\lambda,
\alpha_1^{\vee}))$.\hfill$\diamondsuit$

In fact the proof provides a more general result. Consider the
following diagramme with natural maps 

$$ \xymatrix{ Z_{s_{r+1}} \ar[d] & P_{r+1}/B \ar[d] \\
 Z_{{\bf w}[r+1]} \ar[r]^{f_{r+1}} \ar[d]^{\psi_{r+1}} & G/B
\ar[d]^{\pi_{r+1}} \\
   Z_{{\bf w}[r]} \ar @/^1pc/ [u]^{\sigma_r} \ar[r] \ar[ur]^ {f_r}&
G/P_{r+1} } $$
then we have:

\blem The relative tangent bundle of $\psi_{r+1}$ is $f^*_{r+1}
(L_{\alpha})$, where $L_{\alpha}$ is the line bundle on $G/B$
associated to the character $\alpha$, and $\alpha$ is the simple root
corresponding to $s_{r+1}$. \hfill $\diamondsuit$
\elem

\section{Construction of the Degeneration}

Consider a sequence ${\bf w} = (s_1, s_2, \cdots , s_m)$ of
simple reflections. Let $\cb$ denote $B\times \ba$ and $\cp_j$  denote  
$P_j\times \ba$ for $1\leq j \leq m$. Then  $\cb^{\bf w} = \cb \times _\ba \cb \cdots, \times_\ba \cb= 
B^{\bf w} \times \ba$ and $\cp^{\bf w}=\cp_1 \times_\ba  \cp_2 \cdots \times_\ba \cp_m  =P^{\bf w}\times \ba$.
Both $\cb^{\bf w}$ and $\cp^{\bf w}$ are group schemes over $\ba$.

\bdefe \label{bottaction} Define the following twisted action of
$\cb^{\bf w}$ on $\cp^{\bf w}$  over $\ba$ as follows:

~~~~$ [(p_1, ~ p_2,~ \cdots,~
p_m),~x] \cdot [(b_1, ~ b_2, ~ \cdots,~ b_m),~x] $ \\
$~~~~~~~~~~~=[(p_1~
b_1,~\phi_x(b_1)^{-1} p_2 b_2, ~ \cdots,~\phi_x(b_{m-1})^{-1} p_m b_m)
,~x] $ 

where $\phi_x:B\to B$ is the family of homomorphisms defined in
Lemma~\ref{familyhom}.  \edefe

\blem  The (right) action of $\cb^{\bf w}$ on $\cp^{\bf w}$ over $\ba$
is {\em free}.

\elem

\noindent{\bf Prooof:}
Recall that an action $\sigma:G\times_S X\longrightarrow X $ of a group scheme $G$ over a scheme $S$ 
on a Scheme $X$ over $S$ is said to be free if the map $(\sigma, \pi_2):G\times_S X\longrightarrow 
X\times _S X$ is a closed immersion. [\cite{mumford}, Page 10].

We show the map $\cb^{\bf w} \times_\ba \cp^{\bf w}\longrightarrow 
\cp^{\bf w}\times_\ba \cp^{\bf w}$ is injective.\\
Suppose 
$$((p_1, p_2, \cdots, p_m)\cdot (b_1, b_2, \cdots, b_m), (p_1, p_2, \cdots, p_m))~~~~~~~~~~~~`$$
$$~~~~~~=((p_1^{\prime}, p_2^\prime, \cdots ,p_m^\prime) \cdot (b_1^\prime, b_2^\prime, \cdots, 
b_m^\prime), (p_1^\prime, p_2^\prime, \cdots, p_m^\prime))$$ over a point $x$. 
Then  $p_1=p_1^\prime, ~ p_2=p_2^\prime, ~\cdots ~p_m=p_m^\prime$,
 and (refer Definition \ref{bottaction}) 
$p_1b_1 = p_1 b_1^\prime, ~ \phi_x(b_1)^{-1}p_2b_2= \phi_x(b_1^\prime)^{-1}p_2 b_2^\prime,$ 
$~ \phi_x(b_2)^{-1}p_2b_3= \phi_x(b_2^\prime)^{-1}p_2 b_3^\prime$, $\cdots,$  
$~ \phi_x(b_{m-1})^{-1}p_mb_m= \phi_x(b_{m-1}^\prime)^{-1}p_m b_m^\prime.$ 
Which  successively implies $b_1=b_1^\prime, ~ b_2=b_2^\prime, ~\cdots ~b_m=b_m^\prime$.
Using the fact that the $B$ orbits in $P_i$ are closed embedding, one can show that the 
map $\cp^{\bf w} \times_\ba \cb^{\bf w}\longrightarrow 
\cp^{\bf w}\times_\ba \cp^{\bf w}$ is a closed immersion.
\hfill$\diamondsuit$

\bdefe \label{bsdschemes}
 Since the action of $\cb^{\bf w}$ on $\cp^{\bf w}$ over $\ba$
is {\em free}, the quotien $ {\cz _{\bf w}} =\cp^{\bf w}/ \cb^{\bf w} $ exists
 as an algenraic space over $\ba$ [\cite{Keel}, Theorem 1.1]. Let $\pi: \cz_{\bf w}\lr \ba$ denote the
defining morphism. 
\edefe

\brem \label{induct} Note that the projection $\cp^{{\bf w}[r]}\to
\cp^{w[r-1]}$ is equivariant for the actions of $\cb^{{\bf w}[r]}$ and
$\cb^{w[r-1]}$ on $\cp^{w[r]}$ and $ \cp^{{\bf w}[r-1]}$
respectively. Hence this projection descend to give a morphism of the
quotient spaces $\psi_{r,\ba} :{\cz}_{{\bf w}[r]} \lr {\cz}_{{\bf
w}[r-1]}$. It is  a $P_r/B$ fibration. There is also
a section $\sigma_{r,\ba} :{\cz}_{{\bf w}[r]} \hookrightarrow
{\cz}_{{\bf w}[r+1]}$ to the projection $\psi_{r+1,\ba}$ induced by
the inclusion
 $$
\cp^{{\bf w}[r]} \cong P_1 \times P_2 \times \ldots
\times P_{r}\times \{1\}\times \ba \hookrightarrow P_1 \times P_2
\times \ldots
\times P_{r+1}\times \ba = \cp^{{\bf w}[r+1]} $$
      \erem

Now we describe a scheme structure on the algebraic space ${\cz}_{\bf w}$. 

\bth \label{piprops} The morphism $\pi: {\cz_{\bf w}} \lr \ba$ has the
following properties: 
\begin{enumerate} 
\item $\pi$ is a smooth projective morphism 
\item The fiber over $1$, $\cz_{\bf w}^1:=\pi ^{-1}(1)$ is the BSDH variety $Z_{\bf w}$ and  
$\cz_{\bf w}^x:=\pi ^{-1}(x)$ is isomorphic
to $Z_{\bf w}$ for $x \neq 0$
\item $\cz_{\bf w}^0:=\pi ^{-1}(0)$ is a smooth toric variety.

\end{enumerate}
\eeth

\noindent {\bf Proof} of (1) We construct $\ba$ schemes $\cy_i$  for $i=0, \cdots, m$ and smooth morphisms
$\psi_{i,\ba}:\cy_i/B\lr \cy_{i-1}/\cb$,  and sections $\sigma_{i-1,\ba}:\cy_{i-1}/B\lr \cy_{i}/\cb$
for $i=1, 2, \cdots, m$ inductively.  Set $\cy_0 =\cb$, $\cy_1 = \cp_1$ and $\psi_{1,\ba}$ is the 
defining morphism from 
$\cz_{{\bf w}[1]}=\cp_1/\cb \lr \ba$.  For $i\geq 2$, define 
 \[ \cy_i := \cy_{i-1} \times^\cb \cp_i \] 
where the $\cb$ action is defined by $(y,p)\cdot b =(yb, \phi_x 
(b)^ {-1}p)$ for $y \in \cy_{i-1}$, $p \in \cp_i$, $b\in \cb$, and $x\in \ba$ .  There is 
still an action of $\cb$ (in fact $\cp_i$) on $\cy_i$ coming from right multiplication on 
$\cp_i$ and $\cy_i/B \cong \cz _{{\bf w}[i]}$. Note that  the $\cp_i$ action is free and $\cy _i/\cp_i\cong \cz_{{\bf w}[i-1]}$. 
% $\cy_i/\cb\cong \cz_{{\bf w}[i]}$
In fact, the  scheme $\cy_i$ is a principal $\cp_i$ bundle over $\cz_{{\bf w }[i-1]}$. 

We define $\psi_{i,\ba}$ to be the composite of the following maps
$$ \psi_{i,\ba}: \cz_{{\bf w}[i]}= \cy_i/\cb \lr \cy_i/\cp_i \cong \cy_{i-1}/\cb = \cz_{{\bf w}[i-1]} $$

The map sending $y$ to the class $[y,1]$ from $\cy_{i-1}$ to $\cy_i$ descend to give 
$$\sigma_{i-1,\ba}:\cy_{i-1}/\cb\lr \cy_i/\cb$$

Since $G$ is simply connected group, there exists an irreducible  rank $2$ 
representation $V_i$ of $P_i$ which trivially extends over 
the base $\ba$ to $\cp_i$. This gives rise to a rank $2$ vector bundle ${\mathcal V}_i$ 
on $ \cz_{{\bf w}[i-1]}$  such that the $\cb$ quotient $\cy_i/\cb\cong \cz_{{\bf w}[i]}$ is 
cannonically isomorphic to ${\bf P}({\mathcal V}_i)$.  We get the following diagramme:

$$\xymatrix{
\cz_{{\bf w}[i]} \ar[d]_{\psi_i, \ba}\ar[r]^\cong & {\bf P}({\mathcal V}_i)\ar[ld]\\
\cz_{{\bf w}[i-1]} &\\
}$$

Since each $\psi_{i, \ba}$ is smooth proper morphism as it is isomorphic to the 
projective bundle of a vector bundle, the composition $\pi=\psi_{m,\ba}\circ
\psi_{m-1, \ba} \circ \cdots \circ\psi_{1,\ba}$ is a smooth projective morphism.

\noindent {\bf Proof} of (2): For $x \neq 0$, consider the map $ f_x : P_1 \times
\ldots \times P_m \lr P_1 \times \ldots \times P_m $ given by
$f_x ((p_1, \ldots, p_m)) = (p_1, \xi(x)^{-1} \cdot p_2, \ldots, \xi (x)^{-1} \cdot 
p_m)$, where $\xi(x) \cdot p_i$ is the multiplication in $P_i$.  We show that this is a 
$B^w$ equivarient map: $f_x ((p_1,p_2\ldots p_m)\cdot  (b_1, b_2, \ldots, b_m))$ =
$f_x ((p_1b_1, \phi_x(b_1)^{-1}p_2b_2, \ldots ,\phi_x(b_{m-1})^{-1} p_mb_m))$
= $(p_1b_1, \xi(x)b_1^{-1}\xi(x)^{-1}p_2b_2 , \ldots , \xi(x)b_{m-1}^{-1}\xi(x)^{-1} p_mb_m)$
(Refer Lemma \ref{familyhom} for $\phi_x$)\\
\noindent
$= (p_1b_1, b_1^{-1}\xi(x)p_2b_2 ,\ldots, b_{m-1}^{-1} \xi(x)p_mb_m))$ 
= $f_x((p_1,p_2  , \ldots , p_m))\cdot  (b_1, b_2, \ldots ,b_m)$
This  $B^{\bf w}$ equivarient isomorphism descends to give a well defined 
isomorphism $\overline{f_x} : Z_{\bf w}
\lr \cz_{\bf w}^x$.   Note that the map $f_1$ is the 
identity map.  Which proves the claim $\cz^1_{\bf w}=Z_{\bf w}$.

\noindent {\bf Proof} of (3):

Note that $\cz_{{\bf w}[i]}^0$ can also be viewed as 
\[ \cz_{{\bf w}[i]}^0 = \cy_{i-1}^0 \times ^{B} (P_i/B) \] where $B$ 
acts on $P_i/B$ via its projection to the maximal torus $T$.  
 We observe that the action of the maximal torus $T$ on
$P_i/B$ factors through the action of the multiplicative group on
the projective line $P_i/B$ via the character $\alpha_i$. We denote
this quotient of $T$ by $\ct_i$.  We define $\ct_{{\bf
w}[r]}:=\ct_1\times\cdots \times \ct_r,$ for all $1\leq r \leq m$.
Then one see that the action of $T\times \cdots \times T$ on
$\cz_{{\bf w}[i]}^0$ factors through $\ct_{{\bf w}[i]}$.

For $i\geq 2$, consider the principal $B$-fibration $\cy_{i-1}^0 \lr
\cy_{i-1}^0/B = \cz_{{\bf w}[i-1]}^0$. Let $E_{i}:= \cy_{i-1}^0 
\times^{B}\ct_i \lr \cz_{{\bf w}[i-1]}^0$ be the principal $\ct_{i}$ 
bundle obtained using the associated construction with the quotient 
homomorphism $B\to \ct_i$. Then $E_{i}$ is a $\ct_{{\bf w}[i]}$ 
variety with $E_{i} \lr \cz_{{\bf w}[i-1]}^0$ is a $\ct_{{\bf w}[i-1]}$ 
equivariant map.  Then we have
$$
 \cz_{{\bf w}[i]}^0 = \cy_{i-1}^0 \times ^{B} P_i/B
= (\cy_{i-1}^0 \times ^{B}\ct_i)\times ^{\ct_i} P_i/B
= E_{i} \times^{\ct_i} P_i/B $$

 Now the Theorem follows as $\cz_{{\bf w}[i]}^0$ has a dense open orbit for the
 action of $\ct_{{\bf w}[i-1]}\times \ct_i\cong \ct_{{\bf w}[i]}$. Hence
 $\cz_{{\bf w}[i]}^0$ is a smooth toric variety.
\hfill$\diamondsuit$

%\blem \label{ele} Let $Z$ be a toric variety with dense torus ${\ct}$
%and $V$ be a ${\ct}\times {\ct}'$-variety for a torus ${\ct}\times
%{\ct}'$. Let $V \lr Z$ be a ${\ct}$-equivariant map which defines a
%principal ${\ct}'$-bundle. Let $X$ be a toric variety with ${\ct}'$
%as its dense torus and consider the associated construction: \[ V
%\times ^{{\ct}'} X = Z' \] Then $Z'$ is a toric variety with dense
%torus ${\ct} \times {\ct}'$.  \hfill$\diamondsuit$ \elem

Notice that $P_{r+1}/B$ has two $\ct_{r+1}$-fixed points, one is the
$B$-fixed point called the {\it Schubert point} and the
other called {\it non-Schubert point}. These give rise to two sections
$$\sigma_r^0~,~\sigma_r^1:~\cz_{{\bf w}[r]}^0\to \cz_{{\bf w}[r+1]}^0 $$
The section $\sigma_r^0$ corresponding to the $B$-fixed
point is called a {\it Schubert section} which is the restriction of
the section $\sigma_{r,\ba}$ to the special fibre.  The other
$\ct_{r+1}$ fixed point of $P_{r+1}/B$ gives the other section
$\sigma_r^1$ disjoint from the Schubert section. This section will be
called {\it non-Schubert section}.

The point $\sigma_{0,\ba}(x) \in (\cp_1/\cb)^x\cong \cz^x_{{\bf w}[1]}$ is called the Schubert
point in $\cz^x_{{\bf w}[1]}$
The {\it Schubert point} of  $\cz^x_{{\bf w}[r]}$ is defined inductively
as the image of the Schubert point under the Schubert section, $\sigma_{r-1,\ba}\mid _ 
{\cz^x_{{\bf w}[r-1]}}$.

A {\it Schubert line} in $\cz_{{\bf w}[r]}^0$ is defined to be any
$\ct_{{\bf w}[r]}$-invariant curve containing the Schubert point. More
generally we may call a face to be a {\it Schubert face} if it contains
the Schubert point.

The following Lemma is standard.

\blem \label{section-quotientAG} Suppose $V$ be a rank $2$ vector bundle over a curve $C$ and 
${\bf P}(V)$ be the projective bundle. Then the sections $\sigma : C \lr {\bf P}(V)$ are in one to 
one corresspondence with the line bundle quotients $Q$ of $V$. Moreover the self intersection number 
of $\sigma(C)$ in the surface is given by $deg ~Q - deg ~V$, where $deg$ denote the degree of the 
locally free sheaves.
\hfill $\diamondsuit$
\elem

We also need the following result which can be extracted from the proof of Theorem~\ref{piprops}.  
We state this as a separate Lemma. 

\blem \label {nonsubert-bundle defn}
The section $\sigma_{i-1,\ba}:\cz_{{\bf w}[i-1]} \lr \cz_{{\bf w}[i]}$ provides a line bundle quotient 
${\mathcal V}_i\lr {\mathcal Q}.$ Let ${\mathcal S}_i$ be the kernel (a line bundle).  For each $x\in \ba$ 
let us denote the restrictions of these bundles on $\cz^x_{{\bf w}[i-1]}$ by ${\mathcal V}_i^x$,  
${\mathcal S}_i^x$ and  ${\mathcal Q}_i^x$ respectively. Then 
$$\cz^x_{{\bf w}[i]}\cong {\bf P}({\mathcal V}_i^x)$$
$$\cz^0_{{\bf w}[i]}\cong {\bf P}({\mathcal Q}_i^0 \oplus {\mathcal S}_i^0)$$
with the projection ${\mathcal V}_i^x\lr {\mathcal Q}_i^x$ providing the Schubert 
section $\sigma^0_{i-1}(\cz^x_{{\bf w}[i-1]})$ for all $x\in \ba$. Moreover the non-Schubert section,  
$\sigma^1_{i-1}(\cz^0_{{\bf w}[i-1]})$, is provided by the 
quotient ${\mathcal Q}_i^0 \oplus {\mathcal S}_i^0\lr {\mathcal S}_i^0.$
\hfill $\diamondsuit$

\elem

We can now state a  generalization of Lemma \ref{basic} to the family.

\blem \label{basic-family} Let ${\bf u}= (s_1, s_2)$ be a sequence
of simple reflections. Then the self intersection number of the section
$\sigma_{1,\ba} (\cz_{s_1}^x)$ in the surface $\cz_{\bf u}^x$ is $(\alpha_2, ~ \alpha_1
^{\vee})$, for all $x \in \ba$. 
\elem

\noindent
{\bf Proof:}  The above Lemma \ref{section-quotientAG} and Remark \ref{nonsubert-bundle defn} shows that the self 
intersection number of $\sigma_{1,\ba} (\cz_{s_1})$ does not change in the family.  Over a general fiber 
this surface is isomorphic to a Schubert surface.  For $x\neq 0$ this follows from  
Lemma \ref{basic} and Theorem \ref{piprops}(2).
\hfill $\diamondsuit$

\section{Chow group of 1-cycles}

Given a variety $X$ defined over $k$ we denote the group of $1$-cycles
on $X$ modulo numerical equivalence by ${\bf N_1}(X)$. Let ${\bf
A_1}(X)$ denote the real vector space ${\bf N_1}(X)\otimes \br$.
Similarly let ${\bf N^1}(X)$ be the group of line bundles modulo
numerical equivalence and ${\bf A^1}(X)$ denote the real vector space
${\bf N^1}(X)\otimes \br$.  Both ${\bf A_1}(X)$ and ${\bf A^1}(X)$ are
finite dimensional by a theorem of Neron-Severi. It is also known
that the intersection pairing ${\bf A^1}(X)\otimes_\br {\bf A_1}(X)\to
\br$ is perfect.

Let us denote the rational curve $\cz_{{\bf w}[1]}^x=\cz_{s_1}^x$ by $L_1$.
The fibre of $\psi_{r,x}=\psi_{r, \ba} \mid _{\cz^x_{{\bf w}[r]}}$ over the 
Schubert point is the Schubert line
$L_{r}$. We index (label) all $\ct_{{\bf w}}$-invariant curves in
$\cz_{{\bf w}}^0$, which project to a Schubert point in $\cz^0 _{{\bf w}[j]}$
for some $j$, by $m$-admissible sequences as follows. For $r\geq 2$, 
 let $I=(i_1, \cdots, i_j)$ be an
$r-1$ admissible sequence and $L_I$ be the corresponding labelled curve
in $\cz_{{\bf w}[r-1]}^0$. Then the $\ct_{{\bf w}[r]}$-invariant curve
$\sigma_{r-1}^0(L_I)$ in $\cz_{{\bf w}[r]}^0$ is denoted by the same
symbol $L_I$ and $\sigma_{r-1}^1(L_I)$ is denoted by $L_{Ir}$, where
$Ir$ denote the $r$-admissible sequence $(i_1,\cdots, i_{j}, r)$.

The group ${\bf A_1}(\cz_{\bf w}^x)$ is freely generated by the Schubert
lines $L_1, ~L_2 ~ \cdots L_m$ (cf. \cite{barton}, Lemma 1.1).
Hence for any $m$-admissible sequence $I=(i_1,i_2\cdots, i_r)$ we have
$L_I= \sum _{j=1}^m d_jL_j,$ for some $d_j\in \br$.  In the next
proposition we give an explicit formula to write down the coefficients
$d_j$.

\noindent {\bf Example 1}: In the following pictures $\ct$-invariant
curves are shown for $Z_{{\bf w}[1]}^0$, $\cz_{{\bf w}[2]}^0$ and
$\cz_{{\bf w}[3]}^0$ respectively.

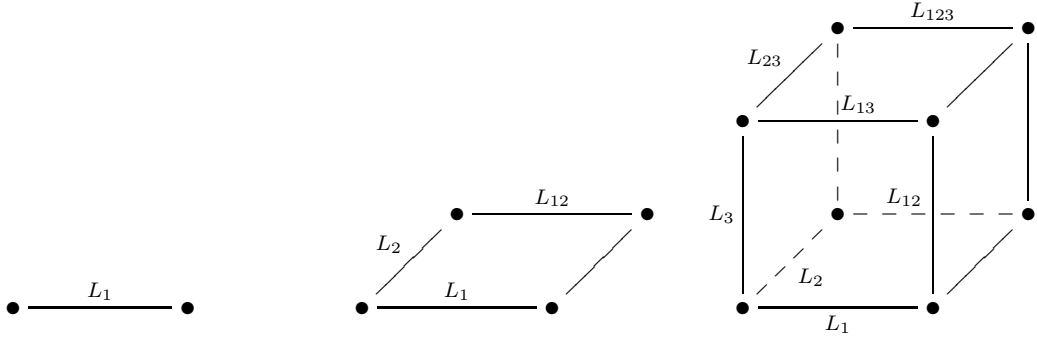
\xymatrix{& & & & & & & & &\bullet \ar@{-}[rr]^{L_{123}}
\ar@{--}[dd] & &\bullet\ar@{-}[dd] \\
& & & & & & & & \bullet
\ar@{-}[ur]^{L_{23}} \ar@{-}[rr]^{~~~~ L_{13}}
\ar@{-}[dd]_{L_3} & &\bullet \ar@{-}[dd] \ar@{-}[ur] & \\
& & & & &\bullet \ar@{-}[rr]^{L_{12}} & &\bullet&  &\bullet
 \ar@{--}[rr]^{\hspace{-.3in}L_{12}} &
& \bullet  \\
\bullet \ar@{-}[rr]^{L_1} & &\bullet & &\bullet \ar@{-}[rr]^{L_1}
\ar@{-}[ur]^{L_2}& &\bullet \ar@{-}[ur] & & \bullet
\ar@{-}[rr]_{L_1}\ar@{--}[ur]_{L_2} &  &\bullet \ar@{-}[ur] & }

\noindent {\bf Example 2}: consider the curve $L_{35679}$ in $\cz_{{\bf
w}[r]}^0$ for $r\geq 9$. This line project down to $L_3$ in $\cz_{{\bf
w}[3]}^0$, i.e
$\psi_4\psi_5\psi_6\psi_7\psi_8\psi_9\cdots\psi_{r}(L_{35679})=L_3$ and \\
$\sigma_{r-1}^0\cdots\sigma_9^0\sigma_8^1\sigma_7^0 \sigma _6^1
\sigma_5^1\sigma_4^1 \sigma_3^0(L_3)=L_{35679}$.

The following observations will be used in later sections:
\begin{enumerate}
 \item The schubert section $\sigma_{j}^0(\cz^0_{{\bf w}[j]})$ and the non-Schubert section
$\sigma_{j}^1(\cz^0_{{\bf w}[j]})$ do not intersect.  In the above picture (Example 1), the 
Schubert surface is placed at the `bottom' and the non-Schubert surface placed at the `top'.

\item Let $I=(i_1, i_2, \cdots, ,i_r)$ be an $m$ admissible sequence.  Then $L_I$ does not exist 
in $\cz^0_{{\bf w}[j]}$ for $j<i_r$,  Moreover $L_I$ lies in the non Schubert section of 
$\cz^0_{{\bf w}[j]}$ for $j=i_r$ and in the Schubert section of $\cz^0_{{\bf w}[j]}$ for all $j>i_r$

\end{enumerate}

\brem (i) The indexing set consists of $2^m-1$ elements as it corresponds
to the set of nonempty subsets of the set $\{1, \cdots , m\}$.

(ii) We show that the total number of $\ct_{\bf w}$-invariant curves
in  $\cz_{{\bf w}}^0$ is $m2^{m-1}$.

One can inductively show that the number of  $\ct_{\bf w[r]}$
invariant points, $pt_r$, in  $\cz_{{\bf w}[r]}^0$ is $2^r$. The number of 
$\ct_{\bf w[r]}$ invariant curves, $l_r$, in  $\cz_{{\bf w}[r]}^0$ can be
counted as follows: There are $l_{r-1}$   $\ct_{\bf w[r]}$-invariant
curves in the Schubert section $\sigma ^0 (\cz_{{\bf w}[r-1]}^0)$
and $l_{r-1}$   $\ct_{\bf w[r]}$-invariant
curves in the non-Schubert section $\sigma ^1 (\cz_{{\bf
    w}[r-1]}^0)$. Also there are $2^{r-1}$   $\ct_{\bf w[r]}$-invariant
curves in the fibre which project down to the $2^{r-1}$  $\ct_{\bf
  w[r-1]}$-invariant points of  $\cz_{{\bf w}[r-1]}^0$.  Thus
$$l_r = 2l_{r-1}+pt_{r-1} ~~ r\geq 2$$

Now we prove the assertion by induction on $m$.   Clearly the
assertion is true when $m=1$.   Assume the result for $m=r$, the
number of $\ct_{\bf w[r]}$-invariant curves in  $\cz_{{\bf w}[r]}^0$ is
$r2^{r-1}$. By the above observation the number of $\ct_{\bf
w[r+1]}$-invariant curves in  $\cz_{{\bf w}[r+1]}^0$ is $2r2^{r-1}+2^r
= 2^r(r+1)$.

We can also count these $\ct_{\bf w}$-invariant curves using our
labellings. We have labelled certain curves with
non-empty-ordered(increasing order) subsets of the set $\{ 1, 2,
\cdots, m\}$.  Clearly there are $2^{m-1}$ curves with the labelling
  set starting with $1$. Any $\ct_{\bf w}$-invariant curves with the
  labelling set starting with $i, ~ i\geq2$ will project down to $L_i,
  $ The curve $L_i$ will project down to a $\ct_{\bf w[i-1]}$-invariant
  point (the Schubert point) of $\cz_{{\bf w}[i-1]}^0$.  The number of
  labelling set starting with $i$ is $2^{m-i}$.  But $\cz_{{\bf
      w}[i-1]}^0$ has $2^{i-1}$ $\ct_{\bf w[i-1]}$-invariant points.
  Any $\ct_{\bf w}$-invariant curves projecting down to any of these
  points other than the Schubert points are not labelled. There are
  $2^{m-i}$ $\ct_{\bf w}$-invariant curves over each of the $\ct_{\bf
    w[i-1]}$-invariant points, each of these curves is equivalent to
  one of those labelled curves.  The total number of curves (labelled and
  un-labelled) which are equivalent to one of the labelled curves with
  labelling set starting with $i$ is $2^{m-i}\times 2^{i-1}=2^{m-1}$.
  This is true for all $1\leq i\leq m$.  Hence the total number
  $\ct_{\bf w}$-invariant  curves is $m2^{m-1}$.

 \erem

\blem  Suppose $Y\subset X$ be smooth projective varieties and $C_1$, $C_2$ be two curves in $Y$.
  If $a_1C_1 +a_2C_2=0$ in ${\bf A}_1(Y)$ then $a_1C_1 +a_2C_2=0$ in ${\bf A}_1(X)$.
\elem

In most of our computations on the limiting toric varieties $\cz_{{\bf w}[r]}^0$ we will work on a 
suitable surface (refer Lemma \ref{basic-family}) $\cz_u^0 \subseteq \cz_{{\bf w}[r]}^0$ . The 
following Lemma is used repeatedly in most of the computation related to curves.

\blem \label{basic-formula} Let ${\bf u}= (s_1, s_2)$ be a sequence
of simple reflections. Let $L_1$ denote the section  $\sigma^0_1(\cz_{s_1}^0)$,  $L_{12}$ denote the 
section $\sigma^1_1(\cz_{s_1}^0)$ and $L_2$ denote the fiber $(P_2/B)^0$ over the Schubert point in $\cz_{\bf u}^0$.
Then $L_{12}=L_1-(2,1)L_2$ in ${\bf N}_1(\cz_{\bf u}^0).$ Recall that the notation $(2,1)$ stands for 
$(\alpha_2, \alpha_1^\vee)$
\elem

\noindent
{\bf Proof:} 
Since the curves $L_1, L_2$ generate ${\bf A}_1(\cz ^0_{{\bf u}})$ there exist numbers $a, b$ such that 
$L_{12} = a L_1 +b L_2$.  Recall that $L_1\cdot L_1 =(2,1)$ (by Lemma \ref{basic-family}),
 and $L_2\cdot L_2=0$ (because $L_2$ is a fiber). Clearly $L_1 \cdot L_{12}=0$, $L_1 \cdot 
L_{2}=1=L_2\cdot L_{12}.$
By intersecting with $L_2 $, $L_{12}\cdot L_2 = a L_1\cdot L_2 +b L_2\cdot L_2,$ 
we get $1=a$.  Similarly intersecting  with $L_1$, we get $0= a(2,1)+b$.  
Thus $a=1$ and $b= -(2,1)$

\hfill $\diamondsuit$

Recall $s_{\alpha^\vee}(x):= x-(x, (\alpha^\vee)^\vee)\alpha^\vee= x-(x,\alpha)\alpha^\vee.  $ 

\bdefe \label{alphacheckdefn}
Set $\alpha_{j_1j_2 \cdots
j_r}^\vee:= s_{\alpha_{j_r}^\vee}s_{\alpha_{j_{r-1}}^\vee}\cdots
s_{\alpha_{j_2}^\vee}(\alpha_{j_1}^\vee)$.
\edefe

\bprop \label{formula} Let ${\bf w}=(s_1, \cdots, s_m)$ be a sequence of
simple reflections with $m\geq 2$. For $2\leq r\leq m$, let
$I=(i_1,i_2,\cdots, i_r):=(i_1,i_2,I')$ be an $m$-admissible sequence
and $L_I$ be the corresponding curve. Then in the Chow group ${\bf
A_1}(\cz_{\bf w}^0)$ of $1$-cycles, we have:

\noindent
(i) $L_{I}=L_{i_1I'}-(i_2, i_1)L_{i_2I'}$ \\
 (ii) $L_I = L_{i_1} +\sum
_{k=2}^r~ \sum_{d=2} ^{k}~ \sum _{1=j_1 < \cdots <j_d =k} (-1)^{d+1}
(i_{j_d}, i_{j_{d-1}})\cdots (i_{j_2}, i_{j_1})L_{i_k}$ \\
(iii) Set $d_{i_1}:=1$ and $d_{i_k}= \sum_{j=1}^{k-1}-d_{i_j}(i_k, i_j)$ for
$2\leq k\leq r$. Then $L_I = \sum_{j=1}^r d_{i_j}L_{i_j}.$ \\
(iv) $d_{i_j}=-(\alpha^\vee_{i_1 i_2 \cdots i_{j-1}}, \alpha_{i_j})$ for all $j>1$
\eprop

\noindent {\bf Proof} {\em of (i):} Note that the curves, $L_{i_1I'}$,
$L_{i_2I'}$ and $L_{i_1i_2I'}$ are edges of a surface that is obtained
from the Schubert surface $<L_{i_1}, L_{i_2}, L_{i_1i_2}>$ by successively
taking sections and hence are isomorphic. Now (i) follows from Lemma
\ref{basic-formula}.

\noindent {\bf Proof} {\em of (ii):}
By (i) $L_{i_1i_2 \cdots i_r}=L_{i_1 i_3\cdots
i_r}-(i_2, i_1)L_{i_2 i_3 \cdots i_r}$.
By induction we get:
$$L_{i_1 i_3\cdots i_r} = $$
$$L_{i_1}- \sum_{k=3}^r (i_k, i_1)L_{i_k} +
\sum_{k=4}^r\sum
_{d=1}^{k-3} \sum _{2<j_1 < \cdots
<j_d <k} (-1)^{d+1}(i_k, i_{j_d})(i_{j_d}, i_{j_{d-1}})\cdots
(i_{j_1}, i_1)L_{i_k}$$

$${\rm and}~ L_{i_2 i_3\cdots i_r} = $$
$$L_{i_2}- \sum_{k=3}^r (i_k, i_2)L_{i_k}
+\sum_{k=4}^r \sum_{d=1}^{k-3}\sum _{2<j_1 < \cdots
<j_d <k} (-1)^{d+1}(i_k, i_{j_d})(i_{j_d}, i_{j_{d-1}})\cdots
(i_{j_1}, i_2)L_{i_k}$$

By multiplying the second equation with $-(i_2,i_1)$ and adding to the
first, one can easily check (ii).

\noindent {\bf Proof} {\em of (iii):} For $k=1, 2$, the formula follows
directly from definition and (ii). For $k\geq 3$,
the coefficient of $L_{i_k}$ in the expansion of $L_I$ is given by (ii):
$$\sum_{d=2}^{k}~ \sum _{1=j_1 < \cdots <j_d =k}
(-1)^{d+1}(i_{j_d}, i_{j_{d-1}})\cdots (i_{j_2}, i_{j_1}) =$$
$$ \sum_{d=2}^{k} -(i_{j_d}, i_{j_{d-1}})  \sum _{1=j_1 < \cdots <j_{d-1}}
(-1)^{d}(i_{j_{d-1}}, i_{j_{d-2}})\cdots (i_{j_2}, i_{j_1})$$
$$= \sum_{d=2}^{k}-(i_{j_d}, i_{j_{d-1}}) d_{i_{j_{d-1}}}$$
That proves (iii).

\noindent {\bf Proof} {\em of (iv):}
Proof is by induction on $j$.  Clearly $d_{i_2}=-(i_2, i_1)=-(\alpha_{i_2}, \alpha^\vee_{i_1})=-(\alpha^\vee_{i_1}, \alpha_{i_2})$. Assume the result for $d_{i_1}, \cdots, d_{i_{j-1}}$. Now \\
$(\alpha^\vee_{i_1 i_2 \cdots i_{j-1}}, \alpha_{i_j})=
( s_{\alpha^\vee_{i_{j-1}}} (\alpha^\vee_{i_1 i_2 \cdots i_{j-2}}), \alpha_{i_j})$  
$=(\alpha^\vee_{i_1 i_2 \cdots i_{j-2}}- (\alpha^\vee_{i_1 i_2 \cdots i_{j-2}}, \alpha_{i_{j-1}})\alpha_{i_{j-1}}^\vee, \alpha_{i_j}) = $ 
$(\alpha^\vee_{i_1 i_2 \cdots i_{j-2}},\alpha_{i_j})- (\alpha^\vee_{i_1 i_2 \cdots i_{j-2}}, \alpha_{i_{j-1}})(\alpha_{i_{j-1}}^\vee, \alpha_{i_j}) =$
$(\alpha^\vee_{i_1 i_2 \cdots i_{j-2}},\alpha_{i_j})+d_{i_{j-1}}(\alpha_{i_{j-1}}^\vee, \alpha_{i_j})$ (by inducation hypothesis) $= (\alpha^\vee_{i_1 i_2 \cdots i_{j-3}},\alpha_{i_j})+d_{i_{j-2}}(\alpha_{i_{j-2}}^\vee, \alpha_{i_j})+d_{i_{j-1}}(\alpha_{i_{j-1}}^\vee, \alpha_{i_j})$
$ = \cdots$ $=(\alpha^\vee_{i_1},\alpha_{i_2})+ d_{i_2}( \alpha_{i_2}^\vee, \alpha_{i_j})+ \cdots d_{i_{j-2}}(\alpha_{i_{j-2}}^\vee, \alpha_{i_j})+d_{i_{j-1}}(\alpha_{i_{j-1}}^\vee, \alpha_{i_j}) = -d_{i_j}$ by (iii).

\hfill$\diamondsuit$
\bcor \label{formula-corollary}
Suppose $I=(i_1,i_2,\cdots, i_r)$ with $\alpha _{i_1}=\alpha_{i_2}$. Then
$$L_I=\sum_{j=2}^r c_{i_j}L_{i_j}$$
a) $c_{i_2} = -1$ and for 
$$k>2, ~ c_{i_k}=  \sum_{d=2} ^{k}~ \sum _{2=j_1 < \cdots <j_d =k} (-1)^{d}(i_{j_d}, i_{j_{d-1}})\cdots (i_{j_2}, i_{j_1})$$ \\

b) $c_{i_j}=(\alpha^\vee_{i_2 i_3 \cdots i_{j-1}}, \alpha_{i_j})$ for all $j>2$

\ecor

\noindent
{\bf Proof:} Replace $i_1$ by $i_2$ and use the fact $(i_2, i_1)=2$ in Proposition \ref{formula} (ii) and (iv).

\brem The coefficients $c_{i_j}$ in corollary \ref{formula-corollary} are negative of the coefficients of $L_{I^\prime}$ where $I^\prime = (i_2, i_3, \cdots , i_r)$

\erem

\brem \label{coefficientusingroots} The coefficient of $L_k$ in $L_{i_1i_2 \cdots i_r}$ vanishes
for all $k\not \in \{i_1, i_2 \cdots i_r\}$. For any
$m$-admissible sequence $(i_1, \cdots i_r, I)$ the coefficient of
$L_{i_j}$ in $L_{i_1 \cdots i_r I}$ and the coefficient of $L_{i_j}$
in $L_{i_1 \cdots i_r }$ are the same for all $j, 1\leq j \leq r$.
\erem

\brem\label{H-surfaces} For an $m$-admissible sequence
$I=(i_1,i_2\cdots,i_r)$, the coefficient $d_{i_j}$ of $L_{i_j}$ in
the expression $L_I=\sum d_{i_j}L_{i_j}$ is the negative of the self
intersection of $\sigma^0_{i_{j}-1}(L_{i_1\cdots i_{j-1}})$ in the
surface $\psi^{-1}_{i_j}(L_{i_1\cdots i_{j-1}})$. From the inductive
definition of the coefficients of $L_I$, the coefficients  may seem
to take arbitrary values but Proposition \ref{formula} (iv) shows that 
these numbers belong to the set $\{(\alpha, \gamma^\vee)\mid
\alpha\in S~ {\rm and}~ \gamma ^\vee $ is any dual root.  In other words,
the coefficients are bounded by $coxeter ~ number$ of the dual root system \erem

\section{A basis for ${\bf A}_1(\cz_{\bf w}^0)$}

The main result of this section is the following basis theorem and the
two algorithms to find out this basis from the given sequence ${\bf w}$.

\bth\label{invariantlines} Let ${\bf w} = (s_1, ~s_2, ~ \cdots, s_m)$
be a sequence of simple reflections.  Then there exist a set of $m$
linearly independent $\ct_{{\bf w}}$ invariant curves, $L_j({\bf w})$
$1\leq j\leq m$, of $\cz_{\bf w}^0$ which generate ${\bf A}_1(\cz_{\bf
w}^0)$ such that every $\ct_{{\bf w}}$ invariant curve lies in the
$\bz_{\geq 0}$ span of this set.  \eeth

\noindent
{\bf Proof:} We choose the generating set inductively. To begin the
induction we note that for $\cz_{{\bf w}[1]}^0\cong \bp^1$, the assertions of
the theorem are valid. Suppose we have chosen a generating set
$L_j({\bf w}[r])$, $1\leq j\leq r$ for $\cz_{{\bf w}[r]}^0$ such that every
$\ct_{{\bf w}[r]}$ invariant curve in $\cz_{{\bf w}[r]}^0$ are non-negatively
generated by $L_j({\bf w}[r])$.  Then a generating set for $\cz_{{\bf w}[r+1]}^0$
can be chosen as

\[ L_j({\bf w}[r+1]) :=
\left \{ \begin{array}{ll} \sigma_r^0(L_j({\bf w}[r])) & if \
\sigma_r^0(L_j({\bf w}[r]))^2 \ \leq 0\ in \  
\psi_{r+1}^{-1}(L_j({\bf w}[r])) \\
\sigma_r^1(L_j({\bf w}[r])) & if \ \sigma_r^0(L_j({\bf w}[r]))^2 \ > 0 \
in \   \psi_{r+1}^{-1}(L_j({\bf w}[r]))

                                \end{array}
                                \right.  \] for $1\leq j\leq r$ and
$L_{r+1}({\bf w}[r+1]):=L_{r+1}$. Here $\sigma_r^0(L_j({\bf
w}[r]))^2$ denote the self intersection number in the surface
$\psi_{r+1}^{-1}(L_j({\bf w}[r]))$.  By induction we know that any
$\ct_{{\bf w}[r]}$ invariant curve in $\cz_{{\bf w}[r]}^0$ is a
positive linear combination of $L_j({\bf w}[r])$.  First observe
from the definition of $L_j({\bf w}[r+1])$, $j=1,\cdots r$, that the
set of curves $\sigma_r^0(L_j({\bf w}[r]))$ and $\sigma_r^1(L_j({\bf
w}[r]))$ are non-negative linear combinations of
 $<L_{r+1}, L_r({\bf w}[r+1]), \cdots   L_1({\bf w}[r+1]) >$.
 Any $\ct_{{\bf w}[r+1]}$-invariant curve in $\cz^0_{{\bf
w}[r+1]}$ is either equivalent to the fibre $L_{r+1}$ or lies in
either of the sections $\sigma_r^i(\cz^0_{{\bf w}[r]})$. Now any
$\ct_{{\bf w}[r+1]}$-invariant curve in either of the sections
$\sigma_r^i(\cz^0_{{\bf w}[r]})$ are non-negative linear combinations
of $\sigma_r^i(L_j({\bf w}[r]))$ by induction hypothesis. Hence they
lie in the $\bz_{\geq 0}$ span of $L_{j}({\bf w}[r+1])$.
\hfill$\diamondsuit$

Note that the curve $L_j({\bf w})$ is represented by an $m$-admissible
sequence with initial entry $j$, i.e., $L_j({\bf w})=L_I$ where
$I=(j,\cdots)$.

We have seen that any $\ct_{{\bf w}}$ invariant curve $L_I$ is a
linear combination of Schubert lines $L_i$ with integer coefficients.
But the coefficients can be negative.  For example if $s_1=s_2$ then
$L_{12} = L_1-2L_2$. The idea is to replace $L_1$ by $L_{12}$ in the
generating set whenever $L_{12}=L_1 +d_2L_2$ and $d_2 <0$. So $L_1 (=L_{12}+2L_2)$ is a positive linear 
combination of $L_{12}$ and  $L_2$ This
observation provides an algorithm to obtain the basis.

\blem \label{compalgo} Let ${\bf w} = (s_1, ~s_2, ~ \cdots, s_m)$ be a
sequence of simple reflections.  \begin{enumerate}

\item If $i_2>1$ is the smallest positive integer such that $s_{i_2} 
=s_1$ then  $L_1({\bf w}[j])$ $=L_1, \forall j, 1\leq j \leq
i_2-1$ and $L_1({\bf w}[i_2])=L_{1i_2}$. If there is no such $i_2$
then $L_1({\bf w})=L_1({\bf w}[k])=L_1$, for all $k\geq 1$.

\item Suppose there exist an $i_2$ such that $s_{i_2} =s_1$. If
$i_3 >i_2$ be the smallest positive integer such that $c_3:=(i_3,
i_2)<0$ then  $L_1({\bf w}[j])=L_{1i_2}, \forall j,  i_2 \leq j
\leq i_3-1$ and $L_1({\bf w}[i_3])=L_{1i_2i_3}$. If there is no such
$i_3$ exists then $L_1({\bf w}) = L_1({\bf w}[k])=L_{1i_2}$, for all
$k\geq i_2$.

\item Let $L_1({\bf w}[i_{r-1}])=L_{1i_2i_3 \cdots i_{r-1}}$ be chosen
inductively for $r\geq 3$. If $i_r>i_{r-1}$ is the smallest positive
integer such that $c_r = (i_r, i_2)-c_3(i_r, i_3)-c_4(i_r,i_4)-
\cdots -c_{r-1}(i_r,i_{r-1})<0$  then  $L_1({\bf
w}[j])=L_1i_2\cdots i_{r-1} , \forall j,
  i_{r-1}\leq j \leq i_r-1$ and
$L_1({\bf w}[r]):=L_{1i_2i_3\cdots i_r}$. If there is no such $i_r$
exists then $L_1({\bf w}) = L_1({\bf w}[k])=L_{1i_2\cdots i_{r-1}}$,
for all $k\geq i_{r-1}$.

\item For $j>1$, we repeat these procedures for the truncated word
$[j-1]{\bf w}$.  In other words $L_j({\bf w})=L_1([j-1]{\bf w})$.
\end{enumerate} \elem

\noindent {\bf Proof:}

We obtain $L_1({\bf w})$ inductively. Clearly $L_1({\bf w}[1])=L_1$
gives the required basis for ${\bf A_1}(Z_{{\bf w}[1]}^0)$. By
Theorem~\ref{invariantlines}  and  Lemma \ref{basic}, $L_1({\bf
w}[2])$ is $L_1$ or $L_{12}$ (note that $\sigma_1^0(L_1)=L_1$ and
$\sigma_1^1(L_1)=L_{12}$) depending on whether $(2,1)$ is non
positive or positive.  From the theory of root system we know that
$(2,1)$ is positive only when $s_2=s_1$ If $s_1 \not = s_j, \forall
j>1$ in the sequence ${\bf w}=(s_1, \cdots ,s_m) $ then $L_1({\bf
w}[k]) =L_1, \forall k \geq 1$

If there exists a $j>1$ such that $s_j=s_1$ then let $i_2>1$
be the smallest positive integer such that $s_1=s_{i_2}.$ As above
we have $L_1({\bf w}[k]) =L_1 \forall k<i_2$ and $L_1({\bf w}[i_2])=L_{1i_2}$.
Again by Theorem~\ref{invariantlines} $L_1({\bf w}[i_2+1])$ is $\sigma_{i_2}^0(L_1({\bf w}[i_2]))$ or
$\sigma_{i_2}^1(L_1({\bf w}[i_2]))$ depending on whether
$\sigma_{i_2}^0(L_1(w[i_2])^2 $  is non positive or positive  in the surface
$\psi^{-1}_{i_2+1}(L_{1i_2}) $.  But
$\sigma_{i_2}^0(L_1(w[i_2])^2 $  is the negative of the coefficient of
$L_{i_2+1}$ in $L_{1i_2(i_2+1)}$ ( i.e if $L_{1i_2(i_2+1)} = L_1+
d_{i_2}L_{i_2}+d_{i_2+1} L_{i_2+1}$ then
$\sigma_{i_2}^0(L_1(w[i_2])^2 $ is $-d_{i_2+1}$). \\
By Proposition~\ref{formula}(iii)  $d_{i_2+1} = -1(1,i_2+1)-d_{i_1}(i_2,
i_2+1)=-(i_2, i_2+1) +2(i_2, i_2+1)=(i_2, i_2+1)$ [as $s_1 =s_{i_2}$
and $d_{i_1}=-2$]\\
 This justifies step 2.  A similar argument will work
for step 3.  Note that $c_r$  is the coefficient of $L_{i_r}$
in $L_{1i_2i_3 \cdots i_{r-1}i_r}$, which is
negative of the self intersection number, $\sigma^0_{i_r-1}(L_{1i_2i_3 \cdots i_{r-1}})^2$, 
in the surface
$\psi_{i_r}^{-1}(L_{1i_2i_3 \cdots (i_r-1)})$.  Now the proof
follows from Theorem~\ref{invariantlines}. \hfill$\diamondsuit$

The above algorithm can be written using root data related to the algebraic group
$G$. This will also prove the Remark \ref{coefficientusingroots}.

Recall also that the height of a positive root $\gamma= \sum_i n_i \gamma_i$,
where $\gamma_i$'s are simple roots, is defined to be the number $\sum
_i n_i$ and is denoted by ht$(\gamma)$ in a root system. 
We will use this definion for the `dual roots'
$\{\alpha ^\vee \}$ in the following lemma. 
Recall $\alpha^\vee = \frac{2\alpha}{(\alpha, \alpha)}$ and the simple reflection 
$s_{\alpha^\vee}(x):= x-(x, (\alpha^\vee)^\vee)\alpha^\vee= x-(x,\alpha)\alpha^\vee.  $ 

\bdefe \label{weylalgo}
Set $\alpha_{j_1j_2 \cdots
j_r}^\vee:= s_{\alpha_{j_r}^\vee}s_{\alpha_{j_{r-1}}^\vee}\cdots
s_{\alpha_{j_2}^\vee}(\alpha_{j_1}^\vee)$.

Let ${\bf w} = (s_1, s_2, \cdots, s_m)$ be a sequence of simple
reflections. We give an algorithm to find the subsequence $I({\bf w})$ that gives the basis element $L_1({\bf w})$

 Let $i_2>1$ be the  smallest positive integer such that $s_{i_2}
=s_1$. Let  $i_3   >i_2$ be the smallest positive integer such that
    ht($\alpha_{i_2i_3}^\vee)>$ht($ \alpha_{i_2}^\vee)$. 

 Suppose we have chosen $i_1, i_{r-1}$
inductively with $r\geq 3$. Let $i_r>i_{r-1}$ be the smallest positive
integer such that  ht(
$\alpha_{i_2i_3i_4\cdots i_r}^\vee)>$ht( $\alpha_{i_2i_3i_4\cdots
  i_{r-1}}^\vee)$. 

\edefe

\bth
Given a sequence ${\bf w} = (s_1, s_2, \cdots, s_m)$ the Schuber lines $L_{I([j-1]{\bf w})}=
L_j({\bf w})$,  $1\leq j\leq m$.

\eeth

\noindent {\bf Proof:} By definition $\alpha^\vee_{i_2 \cdots
i_j}=s_{\alpha_{i_j}^\vee}(\alpha^\vee_{i_2 \cdots
i_{j-1}})=\alpha^\vee_{i_2 \cdots i_{j-1}}-(\alpha^\vee_{i_2 \cdots
i_{j-1}}, ~\alpha_{i_j})\alpha_{i_j}^\vee$.  It is easy to see that
$(\alpha^\vee_{i_2 \cdots i_{j-1}}, ~\alpha_{i_j})
=c_{i_j}$. Note that  ht($\alpha_{i_2i_3i_4\cdots i_j}^\vee)>$ht( $\alpha_{i_2i_3i_4\cdots
  i_{j-1}}^\vee)$ if and only if $c_{i_j}<0$. Now the Lemma follows from the previous 
Lemma~\ref{compalgo}.
\hfill$\diamondsuit$

\brem Even though the length of the sequence could be very large, the
number of steps in the above algorithm can not be very large.  For
example if the height of the longest root in the dual Weyl group orbit
of $\alpha_1^\vee$ (i.e the unique dominant weight in the orbit
containing $\alpha_1^\vee$) is $n_1$ then the number of steps in the
algorithm for $L_1({\bf w})$ will be at most $n_1$. In other words if
$L_1({\bf w})=L_{1i_2\cdots i_r}$, then $r\leq n_1$. \erem

\section{Extremal Rays and Mori Rays}

Let $X$ be a normal projective variety defined over $k$. Then we recall the
following definitions (motivated by the definition given for $k=\bc$
in Page 254 of \cite{wisniewski}).

\bdefe Let $NE(X)\subset A_1(X)$ be the $\br_{\geq 0}$ cone
spanned by effective $1$-cycles.  A ray $R\subset NE(X)$ is an
{\em extremal ray} if given $Z_1,Z_2\in NE(X)$ such that $Z_1+Z_2\in R$,
then both $Z_1,Z_2\in R$.  \edefe

\bdefe  If an extremal ray $R$ satisfies $R\cdot K_X<0$, then $R$ is
called a {\em Mori extremal ray} (also referred as a {\em Mori ray})
where $K_X$ denote the canonical bundle of $X$. \edefe

\blem\label{extremalrays} Let $X$ be a variety such that $NE(X)$ is
nonnegatively generated by a linearly independent set of effective
$1$-cycles. Then the rays defined by this generating set are precisely
the extremal rays.  \elem

\noindent {\bf Proof:} Let $Z_1,\cdots, Z_n$ be the linearly
independent set of effective curves generating $A_1(X)$. First we
prove that the rays $\br _{\geq 0}Z_i$'s are extremal rays.  Suppose
$\sum_{a_i \geq 0} a_iZ_i+\sum_{b_i\geq 0}b_iZ_i=cZ_k\in \br _{\geq
0}Z_k$. Then $a_i=b_i =0,$ for $i\neq k$ as $Z_i$'s are linearly
independent. Now both $\sum_{a_i \geq 0} a_iZ_i = a_kZ_k$ and
$\sum_{b_i \geq 0} b_iZ_i =b_kZ_k$ lie in $\br _{\geq 0}Z_k$.

Consider any extremal ray generated by $R=\sum_{a_i \geq 0} a_iZ_i$.
If it is not one of the $\br _{\geq 0}Z_i$'s, then there are at least
two nonzero coefficients, say $a_k$ and $a_l$.  Let $R_1:= \sum_{i\neq
l}a_iZ_i$ and $R_2:=a_lZ_l$. Then clearly $R_1+R_2 =R$ lies in $\br
_{\geq 0}R$, but neither $R_1$ nor $R_2$ lies in $\br _{\geq 0}R$, a
contradiction.  Which proves that $\br _{\geq 0}Z_i$'s are the only
extremal rays.\hfill$\diamondsuit$

\bth \label{Extrmalrays} The extremal rays of the toric variety
$\cz_{\bf w}^0$ are precisely the curves $L_j({\bf w})$. \eeth

\noindent {\bf Proof:} In view of Lemma~\ref{extremalrays} it suffices
to prove that the effective cone $NE(Z^0_{{\bf w}})$ coincides with the
positive convex cone generated by the torus invariant curves $L_i{({\bf w})}$. Since
$L_i{({\bf w})}$ form a basis for $A_1(Z^0_{{\bf w}})$, we can write any
effective curve $C$ as a linear combination $\sum n_iL_i{({\bf
w})}$. For each $1\leq l\leq l({\bf w})$, consider the divisor
$D(a_l)=\sum_{i\neq l}{\tilde D}_i({\bf w}) + a_l{\tilde D}_l({\bf
w})$. Then $D(a_l)$ is ample for all $a_l>0$ by
Theorem~\ref{amplecone}(a). Hence $D(a_l)\cdot C =\sum_{i\neq l}n_i +
a_ln_l>0$. But $\sum_{i\neq l}n_i +a_ln_l>0$ for all $a_l>0$ implies $n_l
\geq 0$. Hence $C$ belongs to the cone generated by $L_i({{\bf
w}})$. \hfill $\diamondsuit$

Now we give a criterion for an extremal ray to be a Mori ray (cf. Page
254 \cite{wisniewski}). 

Recall the following standard lemma:
\blem \label{intersectionlemma} Let $Z$ be a complex manifold and $X$, $Y$ submanifolds of $Z$ intersecting transversally.
Let ${\mathcal N}_{X/Z}$ denote the normal bundle of $X$ in $Z$.  Then  ${\mathcal N}_{X/Z}\mid
_{X \cap Y} \cong  {\mathcal N}_{X\cap Y/Y}.$ If $X$ is a divisor in $Z$ then ${\mathcal N}_{X/Z} = 
{\mathcal O}_Z(X)\mid_ X.$
\elem

We recall the boundary of a toric variety as the
complement of the dense open orbit or equivalently as the union of all
torus invariant divisors. Then the boundary of  $Z_{\bf w}^0$ can be inductively shown 
to be  $$\partial(Z_{{\bf
w}[r]}^0):=\psi_{r}^{-1}(\partial(Z_{{\bf w}[r-1]}^0))\cup
\sigma_{r-1}^0(Z_{{\bf w}[r-1]}^0)\cup \sigma_{r-1}^1(Z_{{\bf
w}[r-1]}^0) $$ Then the canonical bundle $K_{Z_{{\bf w}[r]}^0}\cong
-\partial(Z_{{\bf w}[r]}^0)$ (cf. \cite{oda}). Given a Schubert line
$L_r$, we compute it's intersection with the boundary components in the following:

\blem \label{intsectionformula}
For $i\neq r-1$, $(\psi^{i+1})^{-1} \sigma_i^0(\cz^0_{{\bf w}[i]}) \cap  \langle L_r, L_{i+1} \rangle 
= L_r$.  Moreover 
\begin{enumerate}
 \item $ i < r-1, (\psi^{i+1})^{-1} \sigma_i^0(\cz^0_{{\bf w}[i]})\cdot L_r = 0 
 = (\psi^{i+1})^{-1}  \sigma_i^1(\cz^0_{{\bf w}[i]})\cdot L_r$
\item $i > r-1,~(\psi^{i+1})^{-1} \sigma_i^0(\cz^0_{{\bf w}[i]})\cdot L_r = 
(i+1, r)~{\rm and}~(\psi^{i+1})^{-1} \sigma_i^1(\cz^0_{{\bf w}[i]})\cdot L_r = 0$
\item $i = r-1, (\psi^{i+1})^{-1} \sigma_i^0(\cz^0_{{\bf w}[i]}) \cdot L_r = 1 =
(\psi^{i+1})^{-1} \sigma_i^1(\cz ^0 _{{\bf w}[i]} )\cdot L_r $  
\end{enumerate}
\elem

\noindent
{\bf Proof:}
Note that $L_{i+1}\not \subseteq (\psi^{i+1})^{-1} \sigma_i^0(\cz^0_{{\bf w}[i]}),$ in fact, 
$L_{i+1}$ is normal to 	$(\psi^{i+1})^{-1} \sigma_i^0(\cz^0_{{\bf w}[i]})$ and all other 
Schubert lines are contained in $(\psi^{i+1})^{-1} \sigma_i^0(\cz^0_{{\bf w}[i]})$
When $i \neq r-1$, $L_{i+1} \neq L_r$.  Hence $L_r $ is in 
$(\psi^{i+1})^{-1} \sigma_i^0(\cz^0_{{\bf w}[i]})$.  One can easily see that  
$(\psi^{i+1})^{-1} \sigma_i^0(\cz^0_{{\bf w}[i]}) \cap \langle L_r, L_{i+1}\rangle =L_r.$
Now by Lemma \ref{intersectionlemma} $(\psi^{i+1})^{-1} \sigma_i^0(\cz^0_{{\bf w}[i]}) \cdot L_r 
$ is the self intersection number $L_r \cdot L_r$ in the Schubert surface $\langle L_r, L_{i+1}\rangle$

When $ i < r-1$, the Schubert surface $\langle L_r, L_{i+1}\rangle$
maps onto $L_{i+1}$ with firber $L_r$.  Hence $L_r. L_r=0$ in this surface.  Which proves 1.
 
When $i > r-1$, the Schubert surface $\langle L_r, L_{i+1}\rangle$
maps onto  $L_r$  with firber  $L_{i+1}$. Now 2) follows from the Lemma \ref{basic-family}.

It is clear that $L_r$ intersect transversally at the Schubert point and at the non Schubert
point of $L_r$ transversally.  This proves 3.

$\hfill{\diamondsuit}$

\bprop \label{CanonicalbundleIntersectionFormula}
$K_{Z_{\bf w}^0}\cdot L_r= -\partial(Z_{{\bf
w}[r]}^0)\cdot L_r =-2 -\sum_{j=r}^{m-1} (j+1,r) $
\eprop

\noindent
{\bf Proof:}
 Note that  $L_r\subset\sigma_{j}^0
(Z_{{\bf w}[j]}^0)$, for each $j\geq r$, we see that
$\sigma_{j}^1(Z_{{\bf w}[j]}^0) \cdot L_r=0$ for $j\geq r$.  Hence we
have 
$$ K_{Z_{\bf w}^0}\cdot L_r = -(\partial(Z_{{\bf w}}^0))\cdot L_r$$
$$= - (
\sigma_{r-1}^0 (Z_{{\bf w}[r-1]}^0)\ + \sigma_{r-1}^1(Z_{{\bf
w}[r-1]}^0))\cdot L_r - \sum_{j=r}^{m-1}({\psi^{j+1}})^{-1}(\sigma_{j}^0(Z_{{\bf w}[j]}^0))\cdot L_r $$ 
The restriction of the pull back divisor $(\psi^{j+1})^{-1}
(\sigma_{j}^0(Z_{{\bf w}[j]}^0))$ to $L_r$ is isomorphic to the normal
bundle of $L_r$ in the surface $<L_r,L_{j+1}>$.  Hence the degree,
$({\psi^{j+1}})^{-1} (\sigma_{j}^0(Z_{{\bf w}[j]}^0))\cdot L_r = (j,r) $ by
Lemma~\ref{basic-family}.  Now $$K_{Z_{\bf w}^0}\cdot L_r= -\partial(Z_{{\bf
w}[r]}^0)\cdot L_r =-2 -\sum_{j=r}^{m-1} (j+1,r) $$

$\hfill{\diamondsuit}$
\blem \label{nonmori}
If an extremal ray in $Z_{{\bf w}[r]}^0$ is not a Mori ray then the
extremal ray lying over this ray cannot be a Mori ray for any
$Z^0_{{\bf w}[j]},~ j>r$. \elem

\noindent {\bf Proof:} By induction it suffices to prove for
$j=r+1$. Let $I$ be an $r+1$-admissible sequence such that $L_I$ is an
extremal ray in $Z_{{\bf w}[r+1]}$. Assume $\psi_{r+1}(L_I)$ is not a
Mori ray. Then we have 
$$K_{Z_{{\bf w}[r+1]}^0} \cdot L_I = K_{Z_{{\bf
w}[r]}^0} \cdot \psi_{r+1}(L_I) - d $$
 where $d$ is the self
intersection of the curve $L_I$ in the surface
$\psi_{r+1}^{-1}(\psi_{r+1}(L_I))$. Since $\psi_{r+1}(L_I)$ is not a
Mori ray, it follows that $K_{Z_{{\bf w}[r]}^0} \cdot
\psi_{r+1}(L_I)\geq 0$. By the construction of extremal rays
(cf. Theorem~\ref{invariantlines}) it follows that $d\leq 0$. Hence
$K_{Z_{{\bf w}[r+1]}}^0 \cdot L_I \geq 0$ as claimed.
\hfill$\diamondsuit$

This leads to a criterion for an extremal ray to be a Mori ray.

\bth \label{mori} An extremal ray $L_I$ is a Mori Ray if and only if
there exists an $r>0$ such that $L_I$ is the Schubert line $L_r$ and
there is at most one $j>r$ such that $(j,r)<0$ and it should be $-1$.
\eeth

\noindent {\bf Proof:} A Schubert line $L_r$ is an extremal ray if and
only if $(j,r)\leq 0$ for all $j> r$ ( refer Lemma \ref{compalgo}, Theorem \ref{Extrmalrays}). $L_r$ is a Mori ray if and
only if it is an extremal ray and $K_{Z_{\bf w}^0}\cdot L_r=
-\sum_{j>r} (j,r) -2 < 0$ i.e., $\sum_{j>r} (j,r) \geq -1$.  So there
is at most one $j>r$ such that $(j,r)<0$ and it must be $-1$.

Let $I=(i_1,\cdots,i_r)$ such that $L_I=\sum d_{i_j}L_{i_j}$ be a
non-Schubert extremal ray and hence $L_I = L_{i_1}({\bf w})$ (by Theorem \ref{Extrmalrays}). 
Notice that $(i_2,i_1)=2$  and
$(j,i_1)\leq 0$ for all $i_1<j<i_2$   by Lemma \ref{compalgo}.   Hence
 $$K_{Z_{{\bf w}[i_2]}^0}\cdot
L_{i_1i_2} = -\sum_{i_1<j<i_2} (j,i_1) -2 + (i_2,i_1) =
-\sum_{i_1<j<i_2} (j,i_2) \geq 0$$
 This shows $L_{i_1,i_2}$ is
not a Mori ray in $Z_{{\bf w}[i_2]}^0$. Now by Lemma~\ref{nonmori} $L_I$ is
not a Mori ray in $Z_{\bf w}^0$. 
$\hfill{\diamondsuit}$

Recall that a smooth projective variety is called {\it Fano} if it's anti-canonical bundle is ample.

\bcor 
For a ginve ${\bf w} = (s_1, s_2,  \cdots, s_m)$  The toric variety $\cz^0 _{\bf w}$ is Fano 
if and only if all Schubert lines are Mori.

\ecor
\noindent
{\bf Prooof:} Suppose all the $L_i$'s are Mori, then by definition of Mori, they are all extremal
and hence by Theorem \ref{Extrmalrays} and Theorem \ref{invariantlines} they generate all the torus
invariant lines positively.  Now the ampleness of $-K_{\cz^0 _{\bf w}}$ follows from Toric Nakai
Criterion. Conversely if $\cz^0 _{\bf w}$ is Fano then $L_j({\bf w})$ (refer \ref{invariantlines})
are Mori. By Theorem \ref{mori} $L_j({\bf w})$  is the Schubert line $L_j$. Thus all Schubert lines 
are Mori.
\hfill $\diamondsuit$

\noindent {\bf Example 3}: $G=SL_{n+1}$ and the sequence ${\bf
w}_0=(1,2,\cdots,n,1,2\cdots,n-1,\cdots.1,2,1)$ is a reduced
expression for the longest element $w_0$ of the Weyl Group. Then
there are exactly $n$ Mori rays which are:\\
 $L_n,
L_{n+n-1},\cdots,L_{n+n-1+n-2+\cdots n-r}, \cdots, L_{n(n+1)/2}$.

\section{Ample Cone}\label{Amplecone}

Let ${\bf w}=(s_1,s_2, \cdots, s_m)$ be a sequence of simple
reflections and $Z_{{\bf w}[i]}$ be the intermediate
BSDH variety. Let $f_i: Z_{{\bf w}[i]} \to
G/B$ be the $B$-equivariant map. Let $\cl(\omega_j)$ denote the line
bundle on $G/B$ corresponding to the fundamental weight $\omega_j$.
%Choose the unique $B$-invariant(Schubert) divisor $D_i$ such that
% ${\mathcal L}(\omega_i)\cong {\mathcal O}_{G/B}(D_i)$ is the line bundle
%on $G/B$ corresponding to the Fundamental weight $\omega_i$.

 Let $\psi^i:=\psi_{i+1}\circ\cdots\circ\psi_m:Z_{\bf w}\to Z_{{\bf
w}[i]}$ be the composite projection. Then define 
$\cl_j:=(\psi^j)^*(f^*_j(\cl(\omega_j))$. Then
Lauritzen and Thomsen \cite{laurit} have proved that $\cl_j$
for $j
=1,\cdots, m$ form a basis for the Picard group of line
bundles, hence a basis for ${\bf A}^1(Z_{\bf w})$. In fact they also
proved that the ample cone is the `strict positive cone' generated
by $\cl_j$.
\blem \label{LTintersection}
\[ \cl _j \cdot L_r = \left \{ \begin{array}{lllcl}
                                0& {\rm if}  & j<r & {\rm or} & \alpha_j \neq \alpha_r\\
                                 1 & {\rm if}  & j\geq r & {\rm and } & \alpha_j = \alpha _r
                              \end{array} \right. \]
\elem
% IF YOU REMOVE THE LAST FULL STOP IN THE ABOVE IT WON'T WORK
\noindent
{\bf Proof:}

For $j<r$, choose a section of $\cl _j$ that does not contain the Schubert point of 
$Z_{{\bf w}[j]}$. Then the inverse image of this section under $\psi^j$ does not intersect with
the Schubert line $L_r$ and hence the $\cl_j\cdot L_r=0$.

For $j\geq r$, $f_m$ restricts to $l_r$ as an embedding wth image is $P_r/B$. Moreover, $\cl_j$ 
restricted to $L_r$ can be identified with the restriction of line bundle  $L(\omega_j)$ to $P_r/B$. 
Now the theorem follows from the fact that $L(\omega_j)$ has degree $0$ if $\alpha_r \neq \alpha_j$
it is $1$ if $\alpha_r = \alpha_j $
\hfill $\diamondsuit$

The boundary components $(\psi_\ba^j)^{-1}\sigma_{{j-1},\ba}(\cz_{{\bf w}[j-1]}),~ 1\leq j \leq m$, are 
simple normal crossing $\ba$ divisors in $\cz _{{\bf w}}$.
 For each $x \in \ba$, these divisors form a basis for the Picard group $Pic(\cz^x_{\bf w}).$  
In particular $\cl _j$ can be expressed uniquely as linear combination of these boundary divisors
of $\cz^1_{\bf w}=Z_{\bf w}.$
$$ \cl _{j}=\sum_{i=1}^j a_{ij}(\psi^i)^{-1}\sigma_{{i-1}}(Z_{{\bf w}[i-1]})$$
Now consider the following relative divisor
$$ \cl _{j, \ba}:=\sum_{i=1}^j a_{ij}(\psi_\ba^i)^{-1}\sigma_{{i-1},\ba}(\cz_{{\bf w}[i-1]})$$

  We denote the line bundle $\cl^0_j$ given by the restriction of this relative divisor to the special 
fiber at $0$.  Now by using Lemma \ref{intsectionformula}, we can show  an analogue of Lemma 
\ref{LTintersection} in the toric variety for the line bundle $\cl^0_j$.

\blem\label{amplenesscriterion}
A relative line bundle   $\cl $  given by 
$\sum_{i=1}^m a_{i}(\psi_\ba^i)^{-1}\sigma_{{i-1},\ba}(\cz_{{\bf w}[i-1]})$
 on $\cz _{\bf w}$ represents a relatively ample bundle
if and only if $\cl^0 \cdot L_j ({\bf w}) >0$ for all $j$, where $\cl^0$ is the line bundle
given by $\sum_{i=1}^m a_{i}(\psi^i)^{-1}\sigma_{{i-1}}(\cz^0_{{\bf w}[i-1]})$
\elem

\noindent {\bf Proof:} Ampleness being an open condition in a flat family to check the relative
ampleness it sufices to check ove the special fiber $\cz ^0_{\bf w}$. Note that by our
choice of $L_j({\bf w})$ any torus invariant curve can expressed as non negative linear 
combination of these  $L_j({\bf w})$ (refer \ref{invariantlines}).  Now the stated condition is 
equivalent to the Toric Nakai Criterion for ampleness for a smooth toric variety.  
\hfill$\diamondsuit$

\bth\label{amplecone}  The ample cone of $\cz _{\bf w}^0$ is a naturally a subcone of the 
ample cone of the  BSDH variety $Z _w$
\eeth

\noindent {\bf Proof:} 
Using the basis $(\psi_\ba^i)^{-1}\sigma_{{i-1},\ba}(\cz_{{\bf w}[i-1]})$ for the relative 
$Pic(\cz _{\bf w})$ and the openness of the ampleness one can identify the ample cone of 
$\cz^0 _{\bf w}$ as a subcone of the ample cone of $\cz^1 _{\bf w}$.

We would like to give more computable comparison of this sub cone in terms of the natural 
generators of the ample cone of the BSDH variety. For this purpose let us denote the coeffecients 
of $L_j({\bf w})$ by $d_{i_j}([j-1]{\bf w})$ for $j=1,\cdots m$. Then we have 
$$L_k({\bf w}) = L_k + \sum_{{i_j} i_j >k }d_{i_j}([k-1]{\bf w}) L_{i_j}$$

Note that $d_{i_j}([k-1]{\bf w}) <0$ for all $i_j$. Now consider an ample  line bundle $\cl = \sum a_i\cl_i$. Then 
$$\cl \cdot L_k({\bf w}) = a_k + \sum a_id_{i_j}([k-1]{\bf w}) $$
where the sum is taken all indices $i,i_j >k$ and $\alpha_i = \alpha_{i_j}$.  We can prove that $a_m>0$ as 
$\cl \cdot L_m({\bf w}) =a_m$.  By (downward)  induction we assume that $a_i>0$ fo all $i>k$. then we observe 
that  $\cl \cdot L_k({\bf w}) = a_k + \sum a_id_{i_j}([k-1]{\bf w}) >0 $ implies 
$$a_k > -\sum a_id_{i_j}([k-1]{\bf w}) >0 $$

This infact proves that the line bundle $\cl$ is ample on the BSDH variety by the Theorem of Lauritzen
and Thomsen. This gives an inclusion of the ample cone. 
\hfill$\diamondsuit$

{\bf Acknowledgement:}  The authors would like to thank the referee for many valuable comments and 
suggestions.

\begin{center} {\bf Appendix}
\end{center}

{\bf 1.} $L_{i_1 i_2 i_3 i_4 i_5}$ \\
$= L_{i_1  i_3 i_4 i_5} -(i_2, i_1)L_{ i_2 i_3 i_4 i_5} $\\
$ =  L_{i_1   i_4 i_5}- (i_3, i_1)L_{  i_3 i_4 i_5}   -(i_2, i_1)L_{ i_2  i_4 i_5}  + (i_3, i_2) (i_2, i_1)  L_{  i_3 i_4 i_5}  $\\
 $=L_{i_1   i_4 i_5}  -(i_2, i_1)L_{ i_2  i_4 i_5} +[- (i_3, i_1)  + (i_3, i_2) (i_2, i_1)]L_{  i_3 i_4 i_5}  $\\
$ =L_{i_1    i_5}-(i_4,i_1) L_{   i_4 i_5}  -(i_2, i_1)[L_{ i_2   i_5} -(i_4, i_2)L_{   i_4 i_5}]$ \\
   $[- (i_3, i_1)  + (i_3, i_2) (i_2, i_1)][L_{  i_3  i_5}- (i_4, i_3)L_{ i_4 i_5}] $\\
$ =L_{i_1i_5}-(i_2, i_1)L_{ i_2   i_5}$+
$[- (i_3, i_1)  + (i_3, i_2) (i_2, i_1)]L_{  i_3  i_5}$+
$[-(i_4,i_1)+(i_4,i_2)(i_2,i_1)+(i_4, i_3)(i_3,i_1)-(i_4,i_3)(i_3, i_2)(i_2,i_1)]L_{i_4i_5}$\\
$=L_{i_1} -(i_5, i_1)L_{i_5}-(i_2, i_1)[L_{i_2}-(i_5, i_2)L_{i_5}]+$\\
$[- (i_3, i_1)  + (i_3, i_2) (i_2, i_1)][L_{  i_3}-(i_5, i_3)L_{  i_5} ]$+\\
$[-(i_4,i_1)+(i_4,i_2)(i_2,i_1)+(i_4, i_3)(i_3,i_1)-(i_4,i_3)(i_3, i_2)(i_2,i_1)][L_{i_4}- (i_5, i_4)L_{i_5}]$\\
$= L_{i_1} -(i_2, i_1)L_{i_2}$+
$[- (i_3, i_1)  + (i_3, i_2) (i_2, i_1)]L_{  i_3} $+
$[-(i_4,i_1)+(i_4,i_2)(i_2,i_1)+(i_4, i_3)(i_3,i_1)-(i_4,i_3)(i_3, i_2)(i_2,i_1)]L_{i_4}$+
$[ -(i_5, i_1)+(i_5, i_2)(i_2, i_1)+ (i_5, i_3)(i_3, i_1)  +((i_5, i_4)(i_4,i_1)- (i_5, i_3) (i_3, i_2)(i_2, i_1)
-(i_5, i_4)(i_4,i_2)(i_2,i_1)-(i_5, i_4)(i_4, i_3)(i_3,i_1)+(i_5, i_4)(i_4,i_3)(i_3, i_2)(i_2,i_1)] L_{i_5}$

Let ${\mathbf L_{i_1 i_2 i_3 i_4 i_5}:= d_{i_1}L_{i_1}+ d_{i_2}L_{i_2}+d_{i_3}L_{i_3}+d_{i_4}L_{i_4}+d_{i_5}L_{i_5}}$
\hfill $\diamondsuit$

{\bf 2.} Suppose $\alpha_{i_1}=\alpha_{i_2}.$ Then we replace $i_1$ with $i_2$ and use $(i_2, i_1)=2 $ in the above expression to get

 $L_{i_1 i_2 i_3 i_4 i_5}= L_{i_2} -2L_{i_2}+$
$[- (i_3, i_2)  + 2(i_3, i_2)]L_{i_3} $+
$[-(i_4,i_2)+2(i_4,i_2)+(i_4, i_3)(i_3,i_2)-2(i_4,i_3)(i_3, i_2)]L_{i_4}$+
$[ -(i_5, i_2)+2(i_5, i_2)+ (i_5, i_3)(i_3, i_2)  +(i_5, i_4)(i_4,i_2)- 2(i_5, i_3) (i_3, i_2)
-2(i_5, i_4)(i_4,i_2)-(i_5, i_4)(i_4, i_3)(i_3,i_2)+2(i_5, i_4)(i_4,i_3)(i_3, i_2) L_{i_5}$

$=  -L_{i_2}$+
$(i_3, i_2)L_{i_3} $+
$[(i_4,i_2)-(i_4,i_3)(i_3, i_2)]L_{i_4}$+
$[ (i_5,i_2)- (i_5, i_3) (i_3, i_2)
-(i_5, i_4)(i_4,i_2)+(i_5, i_4)(i_4,i_3)(i_3, i_2)] L_{i_5}$

Let ${\mathbf L_{i_1 i_2 i_3 i_4 i_5}:=  c_{i_2}L_{i_2}+c_{i_3}L_{i_3}+c_{i_4}L_{i_4}+c_{i_5}L_{i_5}}$
\hfill $\diamondsuit$

{\bf 3.} $\alpha _{i_1}^\vee$:  ~ $(\alpha _{i_1}^\vee , \alpha _{i_2}) = (\alpha _{i_2} ,\alpha _{i_1}^\vee) =(i_2, i_1) = -d_{i_2}$\\

$\alpha _{i_1 i_2}^\vee = s_{\alpha_{i_2}}(\alpha_{i_1}^\vee)= \alpha_{i_1}^\vee -(\alpha_{i_1}^\vee, \alpha_{i_2})\alpha_{i_2}^\vee :$ \\
$(\alpha _{i_1 i_2}^\vee , \alpha _{i_3})=  (i_3, i_1)-(i_3, i_2)(i_2,i_1) = -d_{i_3}$\\

$\alpha _{i_1 i_2 i_3}^\vee= s_{\alpha_{i_3}^\vee}(\alpha _{i_1 i_2}^\vee)=\alpha_{i_1}^\vee -(\alpha_{i_1}^\vee, \alpha_{i_2})\alpha_{i_2}^\vee -(\alpha_{i_1}^\vee -(\alpha_{i_1}^\vee, \alpha_{i_2})\alpha_{i_2}^\vee, \alpha_{i_3})\alpha_{i_3}^\vee $\\
$=\alpha_{i_1}^\vee -(\alpha_{i_1}^\vee, \alpha_{i_2})\alpha_{i_2}^\vee -(\alpha_{i_1}^\vee , \alpha_{i_3})\alpha_{i_3}^\vee   + (\alpha_{i_1}^\vee, \alpha_{i_2})(\alpha_{i_2}^\vee, \alpha_{i_3})\alpha_{i_3}^\vee :$\\

$(\alpha _{i_1 i_2 i_3}^\vee, \alpha_{i_4})=(i_4, i_1)-(i_4, i_2)(i_2, i_1)-(i_4, i_3)(i_3, i_1)+(i_4, i_3)(i_3,i_2)(i_2,i_1)=-d_{i_4}$\\

It is easy to check that 

$(\alpha _{i_1 i_2 i_3 i_4}^\vee, \alpha_5 ) = -d_{i_5}$ and \\

    $(\alpha _{ i_2  }^\vee, \alpha_3 ) = -c_{i_3}$, $(\alpha _{ i_2 i_3 }^\vee, \alpha_4 ) = -c_{i_4}$,  $(\alpha _{ i_2 i_3 i_4}^\vee, \alpha_5 ) = -c_{i_5}$\\

\noindent

\end{document}